\documentclass{sn-jnl__modified}

\algtext*{EndIf}
\algtext*{EndFor}



\usepackage{tabularx}
\usepackage{csquotes}

\usepackage{xspace}
\newcommand\GAP{\texttt{GAP}\xspace}
\newcommand\Perseus{\texttt{Perseus}\xspace}
\newcommand\CHomP{\texttt{CHomP}\xspace}
\newcommand\RedHom{\texttt{RedHom}\xspace}
\newcommand\polymake{\texttt{polymake}\xspace}
\newcommand\PolymakeJL{\texttt{Polymake.jl}\xspace}
\newcommand\regina{\texttt{Regina}\xspace}
\newcommand\cplus{\texttt{C++}\xspace}
\newcommand\perl{\texttt{Perl}\xspace}
\newcommand\julia{\texttt{Julia}\xspace}
\newcommand\Dionysus{\texttt{Dionysus}\xspace}
\newcommand\PHAT{\texttt{PHAT}\xspace}

\newcommand\FF{\mathbb{F}}
\newcommand\RR{\mathbb{R}}

\newcommand\ZZ{\mathbb{Z}}
\newcommand\Sph{\mathbb{S}}
\newcommand\RP{\mathbb{RP}}

\DeclareMathOperator\AK{AK} 
\DeclareMathOperator\sd{sd}
\DeclareMathOperator\DunceHat{D}
\DeclareMathOperator\Sawblade{SB}
\DeclareMathOperator\SawSquare{SQ}

\jyear{2021}

\theoremstyle{thmstyleone}

\newtheorem{thm}{Theorem}

\theoremstyle{thmstylethree}

\newtheorem{deff}[thm]{Definition}

\theoremstyle{thmstyletwo}

\newtheorem{rem}[thm]{Remark}
\newtheorem{example}[thm]{Example}

\usepackage[normalem]{ulem}

\begin{document}

\title[Frontiers of sphere recognition in practice]{Frontiers of sphere recognition in practice}

\author[1,2]{\fnm{Michael} \sur{Joswig}}\email{joswig@math.tu-berlin.de}

\author[1]{\fnm{Davide} \sur{Lofano}}\email{lofano@math.tu-berlin.de}

\author[1]{\fnm{Frank H.} \sur{Lutz}}\email{lutz@math.tu-berlin.de}

\author[3]{\fnm{Mimi} \sur{Tsuruga}}\email{mimi@elastic.co}

\affil[1]{\orgdiv{Institut f\"ur Mathematik}, \orgname{TU Berlin}, \orgaddress{\street{Stra\ss e des 17.~Juni 136}, \city{10623 Berlin}, \country{Germany}}}
\affil[2]{\orgname{Max-Planck-Institut f\"ur Mathematik in den Naturwissenschaften}, \orgaddress{\city{04103 Leipzig}, \country{Germany}}}
\affil[3]{\orgname{Elastic}, \orgaddress{\street{800 W El Camino Real, Suite~350}, \city{Mountain View}, \state{CA~94040}, \country{USA}}}

\abstract{Sphere recognition is known to be undecidable in dimensions five and beyond, and no polynomial time method is known in dimensions three and four.
  Here we report on positive and negative computational results with the goal to explore the limits of sphere recognition from a practical point of view.
  An important ingredient are randomly constructed discrete Morse functions.}

\keywords{PL manifolds; discrete Morse theory; computational topology}

\maketitle

\section{Introduction}
\noindent
To tell whether a given space is homeomorphic to the sphere in a given dimension is a basic problem in computational topology.
However, this is difficult in an essential way.
\begin{thm}[{S.~P.~Novikov \cite{VolodinKuznetsovFomenko1974}; cf.\ \cite{ChernavskyLeksine2006}}]
  \label{thm:undecidable}
  Given a $d$-dimensional finite simplicial complex $K$ it is undecidable to check if $K$ is homeomorphic to $\Sph^d$ for $d\geq 5$.
\end{thm}
We will consider closed manifolds encoded as finite abstract simplicial complexes---but the methods and results in this article also hold for more general cell complexes with little modification.
For a brief historical overview: $d$-sphere recognition is trivial in dimensions $d\leq 2$.
Rubinstein~\cite{Rubinstein95} and Thompson~\cite{Thompson94} proved that $3$-sphere recognition is decidable.
Subsequently, Schleimer~\cite{Schleimer11} showed that $3$-sphere recognition lies in the complexity class NP, and Lackenby \cite{Lackenby:2021} proved that this problem also lies in co-NP; see \cite{Lackenby:2002.02179} for a recent survey.
The complexity status of $4$-sphere recognition is open.
Summing up, we do not know of any efficient algorithm for $d$-sphere recognition in the relevant dimensions $d\geq 3$.

Our point of departure is that the sphere recognition problem does not go away simply because it is algorithmically intractable.
To the contrary it appears naturally, e.g., in the context of \emph{manifold recognition}, which is the task of deciding whether a given simplicial complex triangulates any manifold and finding its type.
In the piecewise linear (PL) category, recognizing whether a given complex triangulates some PL manifold can be reduced to \emph{PL sphere recognition} since the links of all vertices of the given complex need to be PL~spheres, sometimes also called standard PL sphere.
This plays a role, e.g., for enumerating all manifolds with few vertices or facets \cite{BrehmKuehnel1992,Tricensus,Burton2011,CasaliCristofori2015,Lutz1999,SulankeLutz2009}; 
for detecting errors in experimental topological constructions \cite{AdiprasitoBenedettiLutz2014pre,SpreerKuehnel2011,LutzTsuruga2013ext}; 
or for meshing \cite{SaucanAppleboimZeevi08}.

In the absence of a general sphere recognition procedure the next best thing are certificates for sphericity and non-sphericity, respectively.
A discrete Morse function, $\mu$, on a finite $d$-dimensional abstract simplicial complex, $K$, may be encoded as an acyclic partial matching in the Hasse diagram of the partial ordering of the faces of~$K$; cf.\ \cite{Forman1998,Forman2002} and \cite{Chari00}.
The \emph{critical} faces are those unmatched, and $(c_0,c_1,\dots,c_d)$ is the \emph{discrete Morse vector} of $\mu$, where $c_k$ is the number of critical $k$-faces.
We call such a discrete Morse vector \emph{spherical} if $c_0=c_d=1$ and $c_k=0$ otherwise.
The relevance for our topic comes from the following key result.
\begin{thm}[Whitehead~\cite{Whitehead1939}; Forman~\cite{Forman1998,Forman2002}]
  \label{thm:whitehead_forman}
  A combinatorial $d$-man\-i\-fold is a PL $d$-sphere if and only if it admits some subdivision with a spherical discrete Morse vector.
\end{thm}
So we propose a heuristic method for sphere recognition which navigates between Theorems~\ref{thm:undecidable} and \ref{thm:whitehead_forman}.
There are a few more obstacles though.
Adiprasito and Izmestiev~\cite{AdiprasitoIzmestiev} showed that a sufficiently large iterated barycentric subdivision of any PL sphere is polytopal (and thus inherits a spherical discrete Morse function from linear programming).
However, in view of Theorem~\ref{thm:undecidable}, there cannot be any a priori bound on the number of barycentric subdivisions required to attain polytopality.
Second, deciding whether a discrete Morse function with at most a fixed number $k$ of critical cells exists is \mbox{NP-hard}~\cite{JoswigPfetsch2006,LewinerLopesTavares2003a}, intractable in the parameter $k$ \cite{burton2016parameterized}, and not even a polynomial approximation is available~\cite{BauerRathod2019}.
Finally, there are combinatorial $d$-spheres that do not admit any spherical discrete Morse function~\cite{BenedettiLutz2013a,BenedettiZiegler2011}.

This article is the considerably expanded full version of the extended abstract~\cite{JLT}.
It is organized as follows.
As our first main contribution, in Section~\ref{sec:procedure}, we present an implementation of a sphere recognition heuristic procedure in \polymake, and demonstrate its efficiency.
In the \polymake project, \perl and \cplus are used as programming languages; our heuristic is implemented in~\cplus.
It is also available through the new \julia interface layer \PolymakeJL \cite{Polymake.jl}, which supports the current \polymake Version 4.5.
Section~\ref{sec:runtimes} comprises comprehensive computational experiments which show that there are many randomly constructed, even fairly large, simplicial complexes for which deciding sphericity is surprisingly easy; this agrees with previous observations~\cite{AdiprasitoBenedettiLutz2014pre,BenedettiLutz2014}.
Moreover, on such input our new approach proves to be superior to, e.g., the $3$-sphere recognition implemented in \regina \cite{regina}, which is a standard tool in computational topology.
Note that \regina's method is a full decision algorithm, while our heuristic may be inconclusive.
However, we are not aware of any triangulation of $\Sph^3$ which cannot be recognized by our method.
Another experiment comes from a census of $4$-manifolds provided by \regina; here our heuristics recognizes about $99.5\%$ of the input as spheres or non-spheres.
Finally, in Section~\ref{sec:limitations}, we explore the limitations of our method.
One outcome is the construction of a new family of $2$-dimensional cell complexes which are contractible, but not collapsible.
These \emph{saw blade complexes} generalize the Dunce hat, and in our experiments they occur naturally as one source of difficulty for recognizing spheres.
Moreover, our computer experiments show that there is a \enquote{horizon} for discrete Morse computations, along with implications to homology computations and computational topology in general.

\section{A heuristic sphere recognition scheme}
\label{sec:procedure}

\noindent
We describe our procedure for sphere recognition and its implementation in \polymake \cite{polymake}.
This is the specification:
\smallskip

\noindent
\textbf{Input:} A $d$-dimensional (finite abstract simplicial) complex $K$ with $n$ vertices and $m$ facets, where a \emph{facet} is a face that is maximal with respect to inclusion.  

\smallskip

\noindent
\textbf{Output:} \textsc{Yes}, \textsc{No}, or \textsc{Undecided}, depending whether $K$ has been recognized as a (standard) PL $d$-sphere.

\smallskip

Our procedure features five steps, labeled (0) through (4).
Discussing the trivial preprocessing Step (0) in some details allows us to introduce the basic terminology and notation.
The core Steps (1), (2), (3), and (4) below together yield Algorithm~\ref{algo}.
  
\begin{algorithm}[th]
  \caption{Sphere recognition heuristic}
  \label{algo}
  \begin{algorithmic}  
    \vspace{1mm}
    \Require{Hasse diagram of combinatorial $d$-manifold $K$, where $d\ge 3$}
    \vspace{1mm}
    \Ensure{Semi-Decision: Is $K$ PL homeomorphic to $\Sph^d$?} 
        \vspace{1.5mm}
        \State \hspace{-1.5em} \textbf{(1)} \emph{compute homology} \;
    \If{homology not spherical}{ \Return NO } 
  \EndIf
  \vspace{1.5mm}
  \For{$N$ rounds}
    \State \hspace{-3em} \textbf{(2)} \hspace{0.5em} \emph{compute random discrete Morse vector} \;
    \If{discrete Morse vector is spherical}{ \Return YES }  
    \EndIf
  \EndFor
  \vspace{1.5mm}  
  \For{$N'$ rounds}
    \State \hspace{-3em} \textbf{(3)} \hspace{0.5em} \emph{perform random bistellar flip or edge contraction}\;
    \If{boundary of simplex is reached}{ \Return YES }
    \EndIf
  \EndFor
  \vspace{1.5mm}
  \State \hspace{-1.5em} \textbf{(4)} \emph{compute and simplify presentation of fundamental group $\pi_1$} \;
   \If{presentation is found to be trivial and $d\neq 4$}{ \Return YES }
   \EndIf
   \If{presentation is found to be non-trivial}{ \Return NO }
   \EndIf
   \vspace{1.5mm}
  \Statex \Return UNDECIDED \;
  \end{algorithmic}
\end{algorithm}

\subsection*{(0) Preprocessing}

To verify whether $K$ is a PL $d$-sphere, there are three elementary combinatorial checks that are useful to perform first.
These checks are fast; their running time is bounded by a low-degree polynomial in the parameters $d$, $n$ and~$m$.
If one of the checks fails, this will serve as the certificate that $K$ is not a sphere.

More precisely, we first check if $K$ is \emph{pure}, i.e., each facet has exactly $d+1$ vertices.
Second, we check if each ridge is contained in exactly two facets, where a \emph{ridge} is a face of dimension $d{-}1$.
Success in these two tests will assert that $K$ is a \emph{weak pseudo-manifold} (without boundary). 
Note that the $0$-di\-men\-sional sphere~$\Sph^0$ is a weak pseudo-manifold of dimension $d=0$ with two isolated vertices.

Third, for $d\geq 1$, we check if the $1$-skeleton of $K$ is a connected graph.
A connected weak pseudo-manifold $K$ of dimension $d=1$ is a polygon and thus triangulates $\Sph^1$.

The  pureness and the weak pseudo-manifold property of a simplicial complex is inherited by all face links; cf.~\cite[Rem.~8]{BagchiDatta1998}.
A (connected) weak pseudo-manifold is a \emph{pseudo-manifold} if it is \emph{strongly connected}, i.e., if any two of its facets can be joined by a sequence 
of facets for which consecutive facets share a ridge.
In particular, a pseudo-manifold of dimension $d=2$ is a triangulation of a closed surface or of a closed
surface with \emph{pinch points} (having multiple disjoint cycles as vertex links).

A $d$-dimensional pseudo-manifold $K$ is a \emph{combinatorial $d$-manifold} if all vertex links of $K$ are PL homeomorphic to the boundary 
of the $d$-simplex or, equivalently, if for every proper $i$-face (with $0\leq i<d$) of $K$ its link is a PL $(d{-}i{-}1)$-sphere; here the case $i=d-1$ ensures
the weak pseudo-manifold property.
This recursive nesting of PL spheres suggests an inductive check of the face links of $K$ by dimension, starting with $1$-dimensional 
links of $(d-2)$-faces and proceeding up until the $(d-1)$-dimensional links of the vertices. 

A connected $2$-dimensional weak pseudo-manifold $K$ whose vertex links are single cycles is a combinatorial $2$-manifold and triangulates a closed surface.
If, additionally, the Euler characteristic of $K$ equals two, then $K$ is $\Sph^2$.

If one of the checks on the links of the overall complex $K$ fails, then $K$ cannot be a standard PL sphere, if one of the checks is left undecided 
this leaves $K$ undecided.

After this preprocessing and an inductive check of the vertex-links we may assume that our input looks as follows:
\smallskip

\noindent
\textbf{Input (modified):} 
Let $K$ be a $d$-dimensional combinatorial manifold, for $d\geq 3$.

\smallskip

The subsequent four steps form the core of our sphere recognition procedure.

\subsection*{(1) Homology computation}

Computing the simplicial homology modules of a finite simplicial complex is a standard procedure, which is implemented, e.g., in \CHomP~\cite{CHomP}, \RedHom~\cite{RedHom} or \polymake \cite{polymake}.
The homology with field coefficients can be determined via applying Gaussian elimination to the (simplicial) boundary matrices.
For finite fields of prime order or the rationals this can be achieved in polynomial time (in the size of the boundary matrices); cf.~\cite{Edmonds1967}.
Similarly, over the integers, a homology computation can be reduced to computing Smith normal forms; cf.~\cite[Ch.~11]{Munkres1984}.
Kannan and Bachem~\cite{KannanBachem1979} gave the first polynomial time Smith normal form algorithm, employing modular arithmetic; see also \cite{Iliopoulos1989}.

Here we employ integer coefficients throughout.
A necessary condition for $K$ to be a sphere (PL or not) is $H_d(K)\cong\ZZ$, and all other (reduced) homology groups vanish.
In this case we say that $K$ has \emph{spherical homology}.

\begin{figure}[t]
  \centering
  \includegraphics[height=6cm]{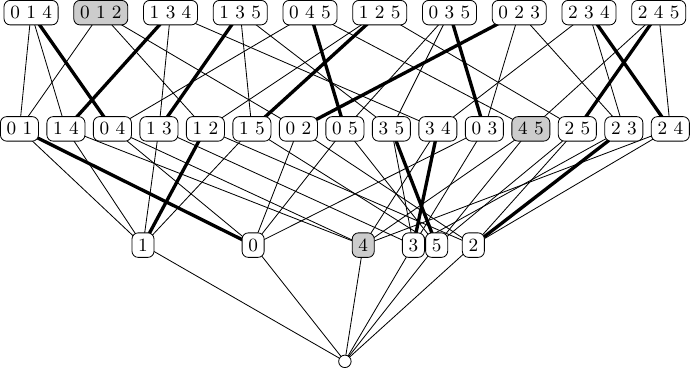}
  \caption{Acyclic matching in the Hasse diagram of $\RP^2_6$. The three unmatched cells are marked.}
  \label{fig:RP2-Hasse}
\end{figure}

The \emph{Hasse diagram} of a simplicial complex $K$ is a directed graph with one node per face of $K$ and a directed arc $(\sigma,\tau)$ if the face $\sigma$ is contained in $\tau$ and $\dim\tau=\dim\sigma+1$.
The non-zero coefficients in the $k$th boundary matrix, which maps $k$-faces to $(k{-}1)$-faces, bijectively correspond to the arcs $(\sigma,\tau)$ in the Hasse diagram for $k=\dim\tau$ (and thus $\dim\sigma=k-1$).
By construction the Hasse diagram is an acyclic graph.

While the modular approach of \cite{KannanBachem1979} and \cite{Iliopoulos1989} is valid for matrices with arbitrary integer coefficients, simplicial boundary matrices have entries $1$, $-1$, and $0$ only.
As a consequence, in an arbitrary simplicial boundary matrix it is \emph{always} possible to perform at least a few Gauss elimination steps.
Moreover, a typical boundary matrix is sparse.
If the matrix happens to stay sparse during the elimination and if, additionally, one does not run out of unit coefficients too soon (such that it is possible to continue with elimination steps), an elimination based Smith normal form algorithm can outperform the more sophisticated modular methods.
This is why for computations of (simplicial) homology elimination algorithms are often preferred; cf.\ Dumas et al.~\cite{Dumas2003} for a survey.

A \emph{(partial) matching} in an arbitrary graph is a subset of the edges such that each node is covered at most once.
In the Hasse diagram of $K$ a matching corresponds to a set of non-zero coefficients in some boundary matrices.
Such a matching, $\mu$, is called \emph{acyclic} if reversing all arcs in $\mu$ (and keeping the arcs not in $\mu$) still gives an acyclic graph.
It is easy to see that an acyclic matching in the Hasse diagram of $K$ yields a sequence of Gauss elimination steps that can be performed in any order without destroying the (unit) pivots required for the subsequent elimination steps.

\begin{figure}[t]
  \centering
  \resizebox{\linewidth}{!}{%
  \begin{tabular}{>{$}c<{$}|*{15}{>{$}c<{$}}}
    \partial_2 & 01 & 02 & 03 & 04 & 05 & 12 & 13 & 14 & 15 & 23 & 24 & 25 & 34 & 35 & \colorbox{gray!20}{45} \\
    \hline
    \mbox{\rule{0pt}{3ex}\colorbox{gray!20}{012}}        &  1 & -1 &    &    &    &  1 &    &    &    &    &    &    &    &    &    \\
    014        &  1 &    &    & \framebox{-1} &    &    &    &  1 &    &    &    &    &    &    &    \\
    023        &    &  \framebox{1} & -1 &    &    &    &    &    &    &  1 &    &    &    &    &    \\
    035        &    &    &  \framebox{1} &    & -1 &    &    &    &    &    &    &    &    &  1 &    \\
    045        &    &    &    &  1 & \framebox{-1} &    &    &    &    &    &    &    &    &    &  1 \\
    125        &    &    &    &    &    &  1 &    &    & \framebox{-1} &    &    &  1 &    &    &    \\
    134        &    &    &    &    &    &    &  1 & \framebox{-1} &    &    &    &    &  1 &    &    \\
    135        &    &    &    &    &    &    &  \framebox{1} &    & -1 &    &    &    &    &  1 &    \\
    234        &    &    &    &    &    &    &    &    &    &  1 & \framebox{-1} &    &  1 &    &    \\
    245        &    &    &    &    &    &    &    &    &    &    &  1 & \framebox{-1} &    &    &  1 \\

    \multicolumn{16}{l}{} \\

    (\partial_1)^{\rm tr}
               & 01 & 02 & 03 & 04 & 05 & 12 & 13 & 14 & 15 & 23 & 24 & 25 & 34 & 35 & \colorbox{gray!20}{45} \\
    \hline
    \mbox{\rule{0pt}{3ex}0}   & \framebox{-1} & -1 & -1 & -1 & -1 &    &    &    &    &    &    &    &    &    &    \\    
    1          &  1 &    &    &    &    & \framebox{-1} & -1 & -1 & -1 &    &    &    &    &    &    \\ 
    2          &    &  1 &    &    &    &  1 &    &    &    & \framebox{-1} & -1 & -1 &    &    &    \\ 
    3          &    &    &  1 &    &    &    &  1 &    &    & -1 &    &    & \framebox{-1} & -1 &    \\ 
    \colorbox{gray!20}{4}          &    &    &    &  1 &    &    &    &  1 &    &    &  1 &    & -1 &    &  1 \\ 
    5          &    &    &    &    &  1 &    &    &    &  1 &    &    &  1 &    &  \framebox{1} &  1 
  \end{tabular}%
  }
  \caption{Boundary matrices of $\RP^2_6$ with coefficients marked that correspond to the acyclic matching of Figure~\ref{fig:RP2-Hasse}; cf.\ Example~\ref{ex:RP2}.}
  \label{fig:RP2-boundary}
\end{figure}

\begin{example}\label{ex:RP2}
  Figure~\ref{fig:RP2-Hasse} shows an acyclic matching, $\mu$, in the Hasse diagram for $K=\RP^2_6$, which is the six-vertex triangulation of the real projective plane.
  Figure~\ref{fig:RP2-boundary} shows the corresponding boundary matrices.
  The pivots corresponding to $\mu$ are marked.
  Using these pivots in an arbitrary order yields an elimination strategy for the computation of the homology modules:
  \[
    \widetilde H_0(\RP^2_6) = 0 \,,\ H_1(\RP^2_6) \cong \ZZ/2\ZZ \,,\ H_2(\RP^2_6) = 0.
  \]
\end{example}

\subsection*{(2) Random discrete Morse functions}
\label{subsec:randomMorse}

A map $f:K\to\RR$ which assigns a real number to each face of $K$ is a \emph{discrete Morse function} if for every $k$-face $\sigma$ of $K$ we have
\begin{equation}\label{eq:discrete-morse}
  \begin{aligned} 
    &\#\{ \tau\in K \, \mid \, f(\tau)\leq f(\sigma),\, \sigma\subset\tau,\, \dim\tau=\dim\sigma+1\} \leq 1 \quad \text{and} \\
    &\#\{\rho\in K \, \mid \, f(\rho)\geq f(\sigma),\, \rho\subset\sigma,\, \dim\rho=\dim\sigma-1\} \leq 1.
  \end{aligned}
\end{equation}
A $k$-face is \emph{critical} with respect to $f$ if both sets in Condition~\eqref{eq:discrete-morse} are empty; the non-critical faces are \emph{regular}, and they form an acyclic matching on the Hasse diagram of $K$; cf.\ \cite{Chari00}.
In this sense the acyclic matchings form equivalence classes of discrete Morse functions.
These concepts were introduced by Forman~\cite{Forman1998,Forman2002}.
The \emph{discrete Morse vector} $(c_0,c_1,\dots,c_d)$ of an acyclic matching counts the critical faces per dimension;
and $K$ is homotopy equivalent to a CW complex with $c_i$ cells in dimension $0\leq i\leq d$.
Let $\FF$ be some field.
A discrete Morse vector is \emph{$\FF$-perfect} for $K$ if $c_i=\beta_i(K;\FF)$ for $0\leq i\leq d$.
\begin{example}
  Figure~\ref{fig:RP2-Hasse} shows an acyclic matching for $\RP^2_6$ with three critical cells.
  The corresponding discrete Morse vector $(1,1,1)$ is $\ZZ_2$-perfect.
  The real projective plane admits a CW-complex structure with one $0$-cell, one $1$-cell, and one $2$-cell.
\end{example}
Discrete Morse vectors which are $\FF$-perfect for any field $\FF$ are \emph{perfect}.
A perfect discrete Morse vector of a sphere, which reads $(1,0,0,\dots,0,1)$, is also called \emph{spherical}.
Theorem~\ref{thm:whitehead_forman} implies that a combinatorial $d$-manifold $K$ that becomes collapsible after the removal of one facet is a PL $d$-sphere.
In 1992, Brehm and K\"uhnel~\cite{BrehmKuehnel1992} used that fact to show that some $8$-dimensional simplicial complex with $15$ vertices is a combinatorial $8$-manifold
(a triangulation of the quarternionic projective plane~\cite{Gorodkov2019}).

By Theorem~\ref{thm:whitehead_forman}, the existence of an acyclic matching whose discrete Morse vector is spherical is a sufficient criterion 
for a combinatorial $d$-manifold $K$ to be a standard PL sphere. This gives rise to the following simple strategy: 
generate discrete Morse functions (or acyclic matchings) at random and check if one of them is spherical; cf.~\cite{BenedettiLutz2014}.

The \texttt{random\_discrete\_morse} function implemented in \polymake has three random strategies which we call \texttt{random-random}, \texttt{random-lex-first}, and \texttt{random-lex-last}.
We will give a short outline and describe the differences among the three strategies and further differences to the original approach from~\cite{BenedettiLutz2014}.

Let $K$ be an arbitrary $d$-dimensional simplicial complex, which is not necessarily a manifold.
A \emph{free face} of $K$ is an $(i{-}1)$-dimensional face that is contained in exactly one $i$-face, $0< i \leq d$.
To save memory, our three strategies are destructive in the sense that they keep changing the complex~$K$.
In each step we pick one of the free faces of codimension one and delete it from $K$ together with the unique $d$-face containing it.
This is an elementary collapse, and the two removed faces form a regular pair, which is a matching edge in the Hasse diagram.
The three strategies only differ in how they pick the free face.
If we run out of free faces, we pick some facet (of maximal dimension), declare it \emph{critical} and remove it.  
After removing a regular pair the dimension of the resulting complex, $K'$, may drop to $d-1$.  
This process continues until $K'$ is zero-dimensional.
In this case, $K'$ only consists of vertices, all of which are declared critical.

For the \texttt{random-random} strategy, we first find all the free faces of $K$ and collect them in a linked list.
If this list is not empty, choose a free face uniformly at random.
Taking the uniform distribution means that each free face has a fair chance of being taken, but this comes at a price since the sampling itself takes time if there are many free faces to choose from.
The reason is that we do not have random access to the free faces, as they are kept in a linked list.
Picking a random element in a linked list takes linear time in the length of that list.
If we run out of free faces, the choice of the critical $d$-face is again uniformly at random.

The strategy \texttt{random-random} is somehow the obvious one, but there is a much cheaper way 
which maintains a certain amount of randomness.  Here the price is that it seems to be difficult to say
something about the resulting probability distribution.  The idea is to randomly relabel the
vertices of~$K$ once, at the beginning, and then to pick the free and critical faces in a
deterministic way (which depends on the random labeling).  Whenever we want to choose a free or critical face,
rather than selecting one at random, we pick the first (in the case of \texttt{random-lex-first}) or
the last one (in the case of \texttt{random-lex-last}) of the linked list.
The \texttt{random-lex-last} strategy was called ``random-revlex'' in~\cite{BenedettiLutz2014}.
We changed the name here to \texttt{random-lex-last} to avoid confusion with the reverse lexicographic (term) ordering, which is different.

The cost of being fair is quite significant, with our current implementation, when dealing with large complexes.
For example, running the \texttt{random-lex-first} and \texttt{random-lex-last} strategies on the fourth barycentric subdivision of $\partial\Delta^4$ took less than three minutes per run whereas the \texttt{random-random} strategy took approximately two hours per run; see Section~\ref{sec:limit_morse_two}.  
It is conceivable that there is some room for improvement here by employing a faster data structure for random sampling; we leave such an implementation for a future version of \polymake.

\begin{rem}\label{rem:homology}
  In Algorithm~\ref{algo}, the Steps (1) and (2) can also be intertwined as finding an acyclic matching results in a partial strategy for computing the homology.
  To this end it is most natural to process the Hasse diagram from top to bottom level by level.
\end{rem}

\subsection*{(3) Random bistellar flips}
\label{subsec:random_flips}
If the previous tests are inconclusive, we can use a local search strategy to determine the PL type; cf.\ \cite{BjoernerLutz2000}.
The boundary $\partial\Delta^{d+1}$ of the $(d+1)$-simplex is a $d$-dimensional simplicial complex with $d+2$ facets.
A \emph{bistellar move} is a local modification of a combinatorial $d$-manifold $K$ in which any subcomplex of $K$ isomorphic to the star of a face in $\partial\Delta^{d+1}$ is replaced by its complementary facets.

To be precise, let $\sigma$ be an $i$-face of $K$ which is contained in exactly \mbox{$d-i+1$} facets $\tau_1,\dots,\tau_{d-i+1}$ 
such that these facets cover exactly $d+2$ vertices. Identifying those $d+2$ vertices with the vertices of $\Delta^{d+1}$ 
yields $(d+2)-(d-i+1)=i+1$ complementary facets $\tau_{d-i+2},\dots,\tau_{d+2}$ in the boundary $\partial\Delta^{d+1}$.
Replacing $\tau_1,\dots,\tau_{d-i+1}$ by  $\tau_{d-i+2},\dots,\tau_{d+2}$ in $K$ is a \emph{candidate bistellar move of dimension $d-i$}, or a \emph{candidate $(d{-}i)$-move} for short.
Let $\sigma' = \cap_{j=d-i+2}^{d+2}\tau_j$ be the complementary face to $\sigma$, where $\sigma'$ is of dimension $d-i$.
If $\sigma'$ is not already contained in $K$, the move is \emph{proper}.
Applying an $i$-dimensional proper bistellar move reduces the $f$-vector of $K$ if and only if $i>d/2$.

Two simplicial complexes are \emph{bistellarly equivalent} if one is obtained from the other by a finite sequence of (proper) bistellar moves.
The following result is essential for the third step in the heuristic.

\begin{thm}[Pachner~\cite{Pachner1987}]\label{thm:pachner}
  A $d$-dimensional simplicial complex is a PL $d$-sphere if and only if it is bistellarly equivalent to $\partial\Delta^{d+1}$.
\end{thm}

This is closely related to Theorem~\ref{thm:whitehead_forman} in the following sense:
Adiprasito and Izmestiev~\cite{AdiprasitoIzmestiev} showed that iterated barycentric subdivisions make any PL sphere polytopal; 
and barycentric subdivisions can be expressed as sequences of stellar subdivisions (which, by Theorem~\ref{thm:pachner}, are connected via sequences of bistellar moves).
Moreover, barycentric subdivisions of polytopal spheres are polytopal, and polytopal spheres admit spherical discrete Morse vectors.

We now discuss the \polymake implementation of the simulated annealing strategy from~\cite{BjoernerLutz2000}.
The function \texttt{bistellar\_simplification} randomly applies bistellar moves to an input of type \texttt{SimplicialComplex} (required to be a combinatorial $d$-manifold) 
with the goal to lower the $f$-vector as much as possible.
In this way the algorithm prefers moves that reduce the $f$-vector; this is called \enquote{cooling}.
It lies in the nature of the sphere recognition problem that we may end up in a local minimum, i.e., when there are no moves to further lower the $f$-vector.
At that point, we deliberately make moves that increase the $f$-vector for some number of rounds (this is called \enquote{heating}).
Then we cool again, hoping that this will help jiggle us out of that local minimum. 
For $0\leq i < \frac{d}{2}$ a corresponding $i$-move will increase the $f$-vector of a triangulation,
while in even dimensions the $f$-vector is not altered by $\frac{d}{2}$-moves.
So in a heating phase we would add vertices via $0$-moves, edges via $1$-moves, etc., and ``randomize'' the triangulation
by performing a (possibly large) number of $\frac{d}{2}$-moves, before returning to cooling. As with all simulated annealing approaches,
adjusting the parameters for the annealing relies on experimentation. For example, we initially may not add vertices in the heating phases
via $0$-moves, as this might successively increase the size of the intermediate complexes, but in case we remain in a local minimum,
first adding some percentage of vertices before performing $1$-moves etc.\ helps in some cases.

Our procedure determines all candidates for bistellar moves of $K$ and sorts them by descending dimension.
During a cooling period we first pick random $d$-moves, if one exists.
Otherwise, we pick random candidate $(d-1)$-moves until we find one which is proper.
If this does not exist either, we continue further to dimensions $d-2$, $d-3$, etc., down to dimension $\lfloor d/2\rfloor+1$.
Any proper move found in this way is applied immediately.
Note that a bistellar move is a local operation, which is why we refrain from copying the entire complex when we apply a bistellar move.
Instead, we perform the operation in place and store the reverse move in a list such that it can be undone later.
Cooling continues until we get stuck with a lexicographically locally minimal $f$-vector.
This ordering of the $f$-vectors is imposed indirectly by preferring higher-dimensional moves.

During a heating period, the story is slightly different. 
One heating strategy is to choose the dimension of the heating move at random with respect to a heuristically determined distribution.
That distribution is encoded as a \emph{heat vector} $(h_0,\dots,h_{\lfloor d/2\rfloor})$ of integers, and we set $h:=h_0+\dots+h_{\lfloor d/2\rfloor}$.
This means, in each round of the heating period we pick the dimension $k$ with probability $h_k/h$, and in that dimension we pick a random proper bistellar move.
For example, the default heat vector in \polymake for $d=4$ is $(10,10,1)$.
This generalizes to the default heat vector $(10,10,\dots,10,1)$ in higher (even) dimensions $d$, while for odd $d$ the pivot dimension $k$ is picked uniformly at random.

Various other parameters control the precise heating behavior; and some of them are adjusted dynamically.
For instance, it is useful to heat up for more rounds if the complex is larger. Sometimes it pays off to experiment with several distributions
or to use other heating schemes. E.g., as mentioned above, first add some (percentage of) or no vertices via $0$-moves, then edges via $1$-moves, etc.

\begin{rem}
\label{rem:bistellar+contractions}
  As a speed-up for large input triangulations, we can first apply edge contractions 
  (with admissible edges for a contraction chosen at random) as long as possible. 
  As experienced for $3$-manifold triangulations~\cite{CsorbaLutz2006}, 
  this eventually leads to a saturation with many edges that block further contractions. 
  Once there is no remaining admissible edge for a contraction,
  we run bistellar flips to reduce the number of edges and then continue with edge contractions again.
  Edge contractions are useful only in an initial phase. Once a local minimum is reached for the size of the triangulation, 
  then bistellar flips are employed to leave the local minimum.
\end{rem}

\subsection*{(4) Fundamental group}

A non-trivial fundamental group $\pi_1(K)$ is a certificate for not being a sphere (PL or otherwise).
Conversely, there is the solution to the (PL) Poincar\'e Conjecture in dimensions other than four.

\begin{thm}[{Smale \cite{Smale1961b}; Perelman \cite{Perelman2002pre}}]\label{thm:PL_poincare}
  Let $K$ be a simply connected combinatorial $d$-manifold, $d\neq 4$, with spherical homology.
  Then $K$ is a PL sphere.
\end{thm}
Freedman proved that a simply connected $4$-manifold with trivial intersection form is homeomorphic to the $4$-sphere~\cite{Freedman1984}.
But his result does not say whether this also holds in the PL category.
In fact, it is a major open problem whether or not \enquote{exotic} $4$-spheres exist.

In~\cite[Chapter~7]{SeifertThrelfall1934} Seifert and Threlfall describe how to obtain a finite presentation of $\pi_1(K)$ from any spanning tree in the $1$-skeleton (with the remaining edges as generators) and all the $2$-faces (as relators).
However, checking if a finitely presented group is trivial is known to be undecidable~\cite{PSNovikov1955}.
Discussing heuristic approaches to simplifying group presentations is beyond the scope of this paper.
In practice we rely on \GAP~\cite{GAP4} which employs Tietze transformations.

Algorithm~\ref{algo} displays our strategy in a concise form; for computational results see Sections~\ref{sec:runtimes} and~\ref{sec:limitations} below.
Notice that the ordering of the Steps (1) through (3) is arbitrary, while the ``YES'' answer in Step~(4) is inconclusive without Step (1).
Yet, there is a benefit from combining Steps (1) and (2); cf.\ Remark~\ref{rem:homology}.
In practice, for a complex $K$ we suspect to be \emph{not} a sphere, we would start with~(1), 
while if we think that $K$ \emph{is} a PL sphere, we first try (2) as a fast routine. If we are 
not successful with (2), we switch to (3), which is slower but can still recognize spheres that do not have
perfect discrete Morse functions; see the discussion in Section~\ref{sec:limitations}.

If, in the case $d\neq 4$, Step (1) gives us a spherical homology vector and Step~(4) a trivial presentation of the fundamental group,
then the overall output is ``YES'', by Theorem~\ref{thm:PL_poincare}. In the case of spherical homology a presentation of the fundamental group 
with only one generator is not possible, but balanced presentations of the trivial group with two generators and two relators
can already be hard to detect; see Section \ref{sec:AkbulutKirby}.

Clearly, when our method gives up with \enquote{UNDECIDED} this does not need to be the end of the story.
For instance, in the $3$-dimensional case we can feed the data into the $3$-sphere recognition procedure of \regina \cite{regina}.
This features a variation of the exact algorithm of Rubinstein \cite{Rubinstein95} and Thompson~\cite{Thompson94}, 
where, for instance, the crushing procedure (the key step of the algorithm) is dramatically simplified \cite{burton2014new}.
In this way \regina can provide certificates for $K$ \emph{not} being spherical based on normal surface theory.
However, we are not aware of a single triangulation of the $3$-sphere for which our procedure fails.

\section{Experiments and runtime comparisons}
\label{sec:runtimes}

\noindent
To find challenging input for Algorithm~\ref{algo} is not entirely trivial.
Most explicit constructions of (standard) PL spheres found in the literature are rather small and can be recognized instantaneously.  
All timings were taken on an AMD Phenom(tm) II X6 1090T Processor CPU (3.2~GHz, 6422~bogomips) and 8~GB RAM with openSUSE Leap 15.0 (Linux 5.1.9-5).

\subsection{Recognizing random $3$-spheres with \polymake and \regina}
A natural class of PL $d$-spheres are the boundaries of $(d+1)$-polytopes obtained as the convex hulls of $n$ points chosen uniformly at random on the unit $d$-sphere in $\RR^{d+1}$.
These have been studied, e.g., in the context of the average case analysis of the simplex method of linear programming \cite{Borgwardt}.
Such examples can be generated with the \texttt{rand\_sphere} command of \polymake.
Table~\ref{tab:rand_sphere} lists \polymake and \regina experiments on $3$-spheres with up to 100,000 vertices.
For more than 15,000 vertices the convex hull computation (necessary only to construct the input) becomes a bottleneck, which is why for the larger examples (marked \enquote{*}) we used connected sums of smaller random spheres.

\begin{table}
\small\centering
\defaultaddspace=0.15em
\caption{Running times (in seconds) on random polytopal $3$-spheres on $n$ vertices.}\label{tab:rand_sphere}
\begin{tabular*}{\linewidth}{@{}r@{\extracolsep{\fill}}r@{\extracolsep{10.50mm}}r@{\extracolsep{10.50mm}}r@{\extracolsep{\fill}}r@{}}
  \toprule
  \addlinespace
  &  \multicolumn{3}{@{}c@{}}{\polymake}      &  \multicolumn{1}{@{}c@{}}{\regina}   \\ \midrule
  \addlinespace
  $n$   &    Morse  &  bistellar &  contr.+bist.    &   \texttt{isThreeSphere}   \\ \midrule
  \addlinespace
    100	&	0.01	&	0.37   &	0.03	&     0.03      \\
    200	&	0.01	&	1.23   &	0.07	&     0.15      \\
    300	&	0.02	&	2.87   &	0.11	&     0.35      \\
    400	&	0.03	&	3.23   &	0.17	&     0.64      \\
    500	&	0.04	&	4.94   &	0.20	&     1.09      \\
    600	&	0.05	&	7.31   &	0.26	&     1.60      \\
    700	&	0.07	&	10.24   &	0.31	&     2.22      \\
    800	&	0.08	&	13.10   &	0.37	&     3.07      \\
    900	&	0.09	&	17.92   &	0.44	&     4.16      \\
   1000	&	0.10	&	23.03   &	0.49	&     5.23      \\
   2000	&	0.38	&       107.85   &	1.25	&     28.22     \\
   3000	&	0.78	&       281.17   &	2.29	&     74.27     \\
   4000	&	1.31	&       551.62   &   	3.41	&     141.65    \\
   5000	&	2.26	&       918.09   &	4.82	&     237.42    \\
  10000	&	8.71	&       4608.71   &	16.48	&    1100.26   \\
  15000	&	22.11	&       /   &		39.77	&    2647.71   \\
 *30000 &	145.90	&       /   &		191.22	&     /      \\
 *50000 &	470.26	&       /   &		515.46	&     /      \\
*100000 &	1586.41	&       /   &		2064.28	&     /      \\
  \addlinespace
  \bottomrule
\end{tabular*}
\end{table}

In all cases, the spheres were successfully recognized by each method.
However, we truncated the time spent on each input to about one CPU hour, such that longer running times are omitted.
The fastest method is \polymake's random search for a spherical discrete Morse function; cf.\ Step~(2) of Algorithm~\ref{algo}.
Nearly competitive is \polymake's procedure of applying edge contractions, combined with random bistellar moves; 
cf.\ Step~(3) of Algorithm~\ref{algo} and Remark~\ref{rem:bistellar+contractions}.

Usually, \regina takes $1$-vertex pseudo-simplicial triangulations as input,  but can also handle (abstract) simplicial complexes.
In the latter case, contracting a spanning tree in the $1$-skeleton yields a $1$-vertex pseudo-simplicial triangulation.
Conversely, the second barycentric subdivision of a pseudo-simplicial complex is a simplicial complex.
In this sense these two encodings of combinatorial manifolds are similar.

\regina's recognition algorithm \texttt{isThreeSphere} runs, as a preprocessing step, the program \texttt{IntelligentSimplify} that transforms the complex into a $1$-vertex 
triangulation and uses bistellar moves to further reduce it, similar to Step (3) of our Algorithm~\ref{algo}. Afterwards the $3$-sphere recognition procedure is employed.
In this way, \regina is able to also find certificates for non-sphericity---which \polymake is incapable of, beyond checking the homology. 
We should also point out that \texttt{IntelligentSimplify} is a heuristic designed to be an out-of-the-box first attempt to simplify a triangulation 
with a polynomial running time, and \regina's bistellar move interface is meant to be interactive. 
This means that with a bit of work to build a custom-made simplification routine, the times in \regina could probably be improved.

The largest simplicial complex in Table~\ref{tab:rand_sphere}, with 100,000 vertices, has 673,274 tetrahedra.
The largest one successfully handled within one hour by \regina has 15,000 vertices and 101,088 tetrahedra. 

Each row of the Tables \ref{tab:rand_sphere}, \ref{tbl:comparison_4d} and~\ref{tbl:comparison_5d6d} corresponds to a single instance only.
However, it is known that there is little variation of, e.g., the number of facets of the convex hull of random points on the unit sphere; cf.\ Reitzner~\cite[Sec.~4]{Reitzner:2005}.
This can also be observed experimentally; cf.\ \cite[Sec.~3.5 and Fig.~6]{polymake:2017} for a closely related setup.
Note that, in fixed dimension $d$, the expected number of facets of a random simplicial $(d+1)$-polytope depends linearly on the number of vertices \cite{Borgwardt}.

\begin{table}[t]
\small\centering
\defaultaddspace=0.15em
\caption{Census of 4-manifolds with up to six pentachora.}\label{tbl:4dim6sim}
\begin{tabular*}{\linewidth}{@{}r@{\extracolsep{\fill}}r@{\extracolsep{\fill}}r@{\extracolsep{5mm}}r@{\extracolsep{\fill}}r@{\extracolsep{5mm}}r@{\extracolsep{\fill}}r@{}}
  \toprule
  \addlinespace
       &  \multicolumn{2}{@{}c@{}}{two pentachora} & \multicolumn{2}{@{}c@{}}{four pentachora} & \multicolumn{2}{@{}c@{}}{six pentachora}  \\ \midrule
  \addlinespace
          &   \# sign.  &  percentage &  \# sign.   &  percentage  &  \# sign.   &  percentage \\ 
  \midrule
 \addlinespace
 Total:  &  8 & $100.0\%$ & 784 & $100.0\%$ & 440,494 & $100.0\%$  \\
 Spheres: &  6 & $75.0\%$ & 642 & $81.9\%$ &  403,240    &  $91.5\%$  \\
 Non-spheres: & 2 & $25.0\%$ & 137 & $17.5\%$ &   35,305     &  $8.0\%$  \\
 Unknown: &  0 & $0.0\%$ & 5 & $0.6\%$ &  1,949      &  $0.5\%$  \\
 \addlinespace
\bottomrule
\end{tabular*}
\end{table}

\subsection{Processing a census of 4-manifolds}
\label{sec:4-manifolds}
We ran our heuristics on a census of $4$-manifolds provided by \regina  \cite{regina}.
These $4$-manifolds are encoded as pseudo-simplicial complexes comprising up to six maximal cells; the $4$-simplices are called \emph{pentachora} 
in \regina for they have five facets (and five vertices).
In combined form, the \texttt{tricensus} command of \regina generates possible facet pairings, then for each such pairing determines all possible gluing permutations, and (on the fly) reduces all gluings to isomorphism signatures that uniquely encode triangulations up to combinatorial isomorphism \cite{Tricensus}.
For the examples with two, four, and six pentachora there are 3, 26, and 639 facet pairings.
And this yields 8, 784, and 440,494 resulting combinatorial types, respectively.
Note that there are no facet pairings for an odd number of maximal cells in even dimensions.
We let \regina expand each of these into a (proper abstract) simplicial complex via the second barycentric subdivision and pass it on to \polymake.
The simplicial complexes resulting from six pentachora have around 4,600 vertices and 86,400 facets. 

Using our heuristic we found the results in Table \ref{tbl:4dim6sim}.
Each positive or negative certificate was obtained in less than four minutes and in 90 seconds on the average.
In all the cases the positive certificates arise from discrete Morse functions, while the negative certificates are provided by non-spherical homology.
Taking row sums in Table \ref{tbl:4dim6sim} we summarize our findings as follows.
\begin{thm}
  Among the $441,286$ combinatorial types of combinatorial $4$-manifolds arising from up to six pentachora\, $91.5\%$ are spheres, and\, $8.0\%$ are non-spheres.
\end{thm}
Thus, our success rate is $99.5\%$, with our heuristic failing on only 1,954 of these combinatorial $4$-manifolds.
All of these have spherical homology; to determine whether they are standard PL $4$-spheres or proper combinatorial homology $4$-spheres 
(combinatorial $4$-manifolds with spherical homology, but not PL homeomorphic to the standard PL $4$-sphere) 
is an interesting question, yet beyond the scope of this article.
Our classification, as a list of \regina's isomorphism signatures of the complexes can be found at \cite{Census}.

\subsection{Higher-dimensional random spheres}
\polymake can easily recognize random polytopal spheres with up to 10,000 vertices in dimension four, 1,000 vertices in dimension five, and 500 vertices in dimension six; cf.\ Tables~\ref{tbl:comparison_4d} and~\ref{tbl:comparison_5d6d}.
Again the input is constructed via uniform random sampling on the unit sphere and taking convex hulls.

\regina provides no heuristic for sphere recognition in dimension four or beyond.
Yet, \regina can simplify a given triangulation of a $4$-dimensional combinatorial manifold via contractions and bistellar moves, returning a smaller pseudo-simplicial complex.
It is not immediate how to check for sphericity from that output.
That implementation is deterministic; thus in each run on a fixed input it gives the same output.
The running times are given in the penultimate column of Table~\ref{tbl:comparison_4d}.
The last column contains the number of simplices remaining after simplification.

\begin{table}
\small\centering
\defaultaddspace=0.15em
\caption{Running times (in seconds) on random polytopal $4$-spheres on $n$ vertices.}\label{tbl:comparison_4d}
\begin{tabular*}{\linewidth}{@{}r@{\extracolsep{\fill}}r@{\extracolsep{5.25mm}}r@{\extracolsep{5.25mm}}r@{\extracolsep{\fill}}r@{\extracolsep{5.25mm}}r@{}}
  \toprule
  \addlinespace
       &  \multicolumn{3}{@{}c@{}}{\polymake} & \multicolumn{2}{@{}c@{}}{\regina}   \\ \midrule
  \addlinespace
  $n$   &    Morse  &  bistellar &  contr.+bist.   &   Simplify  &  number of facets \\ 
  \midrule
  \addlinespace
    100	&	0.04	&	8.22      &	0.46	&     1.52    &		26   \\
    200	&	0.13	&	33.50     &	1.28	&     8.44    &		8    \\
    300	&	0.30	&	76.62     &	2.63	&     24.97   &		42   \\
    400	&	0.54	&	136.85    &	4.77	&     54.51   &	 	78   \\
    500	&	0.82	&	224.67    &	6.17	&     92.25    &	60   \\
    600	&	1.21	&	418.50    &	8.06	&     121.24   &	120   \\
    700	&	1.64	&	639.45    &	10.95	&     184.11   &	98   \\
    800	&	2.28	&	842.94    &	15.38	&     303.16   &	180   \\
    900	&	2.88	&	1109.43   &	16.74	&     370.85   &	144   \\
   1000	&	3.51	&	1418.25   &	22.20	&     474.72   &	170   \\
   2000	&	10.86	&	/         &	40.93	&     2427.81   &	562   \\
   3000	&	26.40	&	/         &	219.44	&     /   &	/   \\
   5000	&	121.90	&	/         &	714.92	&     /   &	/   \\
   10000&	594.46	&	/         &	2633.70	&     /   &	/   \\
  \addlinespace
  \bottomrule
\end{tabular*}
\end{table}

\begin{table}[t]
\small\centering
\defaultaddspace=0.15em
\caption{Running times (in seconds) on random polytopal $5$- and $6$-spheres.}\label{tbl:comparison_5d6d}
\begin{tabular*}{\linewidth}{@{}r@{\extracolsep{\fill}}r@{\extracolsep{6mm}}r@{\extracolsep{\fill}}r@{\extracolsep{6mm}}r@{}}
  \toprule
  \addlinespace
  &  \multicolumn{2}{@{}c@{}}{\polymake, $d=5$}  &  \multicolumn{2}{@{}c@{}}{\polymake, $d=6$}   \\ \midrule
  \addlinespace
  $n$   &    Morse  &  contr.+bist.   &   Morse  &  contr.+bist.  \\
  \midrule
  \addlinespace
    100    &     0.33  	&    10.29       &  10.01   &   301.62     \\
    200    &     1.86 	&    40.52       &  86.78   &   2634.24     \\
    300    &     5.78 	&    102.32      &  387.90  &   /     \\
    400    &     11.49	&    169.86      &  967.55   &  /  \\
    500    &     21.95  &    340.08      &  1788.44   &   /      \\
    600    &     35.31  &    515.86      &    /   &     /    \\
    700    &     55.55  &    820.34      &    /  &      /   \\
    800    &     78.48  &    1120.07     &    /   &     /    \\
    900    &     104.28 &    1441.08     &    /   &     /    \\
   1000    &     133.34 &    2016.53     &    /   &     /    \\
  \addlinespace
  \bottomrule
\end{tabular*}
\end{table}

\subsection{A collapsible $5$-manifold which is not a ball}

We consider the $5$-dimensional simplicial complex $C$ with face vector $f(C)=(5013,72300,290944,495912,383136,110880)$ 
constructed in \cite[Sec.~4]{AdiprasitoBenedettiLutz2014pre}; there $C$ is called \texttt{contractible\_non\_5\_ball}.
This is the first explicit example of a non-PL triangulation of a collapsible (and thus contractible) $5$-manifold, other than the $5$-ball.
By construction, $C$ is a manifold with boundary.
To check the remaining topological properties computationally poses an interesting challenge.

First, the perfect Morse vector $(1,0,0,0,0,0)$ for $C$ was originally obtained in a single random discrete Morse vector search over 82 hours with a \GAP implementation.
The current implementation in \polymake produces the same result (in most runs) in only 9 seconds with the \texttt{random-lex-first} and \texttt{random-lex-last} strategies and in about 10 minutes with the \texttt{random-random} strategy.
This certifies that $C$ is collapsible.

Second, the boundary complex $\partial C$ with face vector $f(\partial C)=(5010,65520,212000,252480,100992)$ was investigated;
it was called \texttt{contractible\_non\_5\_ball\_boundary} in \cite{AdiprasitoBenedettiLutz2014pre}.
Checking all face links for spherical discrete Morse vectors confirmed that $\partial C$ is a combinatorial $4$-manifold.
For each face link a single random try sufficed. In total, the recognition of all face links took about 7.5 hours.
Checking the homology reveals that $\partial C$ is a homology $4$-sphere.
Finally, \GAP identifies the fundamental group $\pi_1(\partial C)$ as the binary icosahedral group.

\section{Limitations}
\label{sec:limitations}

\noindent
In the previous section we saw that many, even fairly large, simplicial spheres can be recognized easily, despite Theorem~\ref{thm:undecidable}.
Here we explore the limitations of our heuristic. The combination of our (positive and negative) experiments may serve as a description 
of a \enquote{horizon} within which we can hope for effective recognition results.

\subsection{General remarks}

We refrain from a detailed comparison of simplicial homology computations. However, standard implementations, such as 
\CHomP~\cite{CHomP}, \RedHom~\cite{RedHom}, \Perseus~\cite{perseus}, and \polymake~\cite{polymake}, 
employ elimination schemes for computing the integer homology, which are equivalent to finding discrete Morse functions with few critical cells.
In this sense, the horizon within which we can compute the simplicial homology is essentially the same as the horizon for 
the discrete Morse Step~(2).
There are more software systems to compute simplicial homology, but many, including, e.g., \Dionysus~\cite{dionysus} and \PHAT~\cite{PHAT}, are restricted to $\ZZ_2$-coefficients.
Finding an optimal discrete Morse function is NP-hard; cf.\ \cite{JoswigPfetsch2006,LewinerLopesTavares2003a}.
Recently Bauer and Rathod established that we may not even hope for polynomial approximability \cite{BauerRathod2019}.

\begin{table}[t]
\small\centering
\defaultaddspace=0.15em
\caption{Collapsing the $d$-simplex.}\label{tbl:dsimplex}
\begin{tabular*}{\linewidth}{@{\extracolsep{\fill}}r@{\hspace{1mm}}r@{\hspace{1mm}}r@{\hspace{1mm}}r@{}}
\toprule
 \addlinespace
 $d$  &  Rounds  &  Non-perfect  &  Percentage  \\ \midrule
   7   &  $10^{10}$  &     0   &   0.0\% \\ 
   8  &  $10^9$  &    12  &   0.0000012\%       \\
   9  &  $10^8$  &     2  &   0.000002\%        \\ 
  10  &  $10^7$  &     3  &   0.00003\%         \\ 
  11  &  $10^7$  &    12  &   0.00012\%         \\ 
  12  &  $10^6$  &     4  &   0.0004\%          \\ 
  13  &  $10^6$  &     6  &   0.0006\%          \\ 
  14  &  $10^5$  &     4  &   0.004\%           \\ 
  15  &  $10^5$  &     8  &   0.008\%           \\ 
  16  &  $10^4$  &     4  &   0.04\%            \\ 
  17  &  $10^4$  &    10  &   0.10\%             \\ 
  18  &  $10^3$  &     2  &   0.2\%             \\ 
  19  &  $10^3$  &     6  &   0.6\%             \\ 
  20  &  $10^3$  &    13  &   1.3\%             \\
  21  &  $10^3$  &    62  &   6.2\%             \\
  22  &  $10^3$  &   153  &   15.3\%             \\
  23  &  $10^2$  &   35  &   35\%             \\
  24  &  $10^2$  &   67  &   67\%             \\
  25  &  $5\cdot 10^1$  &   46  &   92\%             \\
 \addlinespace
\bottomrule
\end{tabular*}
\end{table}

In the subsequent we will exhibit several scenarios in which finding a spherical discrete Morse function for a simplicial sphere may fail in practice.
An obvious impediment is the lack of any spherical discrete Morse function.
The smallest known example is an $18$-vertex triangulation of $\Sph^3$, constructed from a triple trefoil knot supported on three edges~\cite{BenedettiLutz2013a}.

In dimension three, the known recognition \emph{algorithms} for the $3$-sphere make use of normal surface theory.
As a consequence of Theorem~\ref{thm:PL_poincare}, a trivial fundamental group suffices to show that a $3$-manifold is, in fact, the $3$-sphere.

\subsection{Akbulut--Kirby spheres}
\label{sec:AkbulutKirby}

A family of standard PL $\Sph^4$-triangulations, $\AK(r)$, for $r\geq 3$, of the \emph{Akbulut--Kirby spheres}~\cite{AkbulutKirby1985} has been constructed in~\cite{LutzTsuruga2013ext}.
In fact, Akbulut and Kirby~\cite{AkbulutKirby1985} gave handlebody decompositions of a family of $4$-manifolds, in the hope of obtaining exotic $4$-spheres.
Yet, later Gompf~\cite{Gompf} and Akbulut~\cite{Akbulut2010} showed that these manifolds are PL homeomorphic to the standard $4$-sphere $\Sph^4$.
We have
\[ f(\AK(r)) \ = \ (176 + 64 r,\, 2390 + 1120 r,\, 7820 + 3840 r,\, 9340 + 4640 r,\, 3736 + 1856 r). \] 
The triangulated Akbulut--Kirby spheres $\AK(r)$ so far constitute the single explicit family of simplicial spheres that we could not recognize 
easily by our heuristic.
More precisely, Step~(2) failed on all complexes $\AK(r)$ for all $r\geq 3$.
Step~(3) worked for $r=3$, but failed for $r\geq 4$.
Steps~(2) and (3) are particularly relevant in dimension four---in all other dimensions, as a consequence of Theorem~\ref{thm:PL_poincare},  they can conceptually be replaced by the combination of Steps (1) and~(4).
We do not know if the spheres $\AK(r)$ admit perfect discrete Morse vectors.
E.g., the smallest one we found for $r=5$ is $(1,2,4,2,1)$, possibly reflecting that we used two winded up $1$-handles in the construction.

What makes the PL $4$-spheres $\AK(r)$ difficult to recognize is that the original handlebody decomposition of Akbulut and Kirby~\cite{AkbulutKirby1985} is based on the non-trivial presentation
\[
  G(r) \ = \ \langle\, x,y\,\mid\, xyx = yxy;\, x^r = y^{r-1}\,\rangle
\]
of the trivial group built in as the fundamental group, i.e., $\pi_1(\Sph^4)$, of $\AK(r)$.
In 100 out of 450 runs we found \GAP to be able to identify $\pi_1(\AK(4))$ as the trivial group.
However, in dimension four (and knowing that the homology is spherical) this only shows that the input is a topological $4$-sphere, without yielding any information on the PL type.
In one run with $r=5$ we obtained the initial presentations $G(r)$ as the \emph{output} of \GAP's simplification procedure.
For $r\geq 4$, none of the Steps (1)--(4) was conclusive to determine that $\AK(r)$ is the standard PL 4-sphere $\Sph^4$.

\subsection{Contractible but non-collapsible subcomplexes} 

A simplicial sphere $K$ admits a perfect discrete Morse function (respectively vector) if and only if there is a facet $\sigma$ of $K$ such that $K-\sigma$ is a collapsible ball.
In this way a key difficulty in finding a perfect discrete Morse vector for $K-\sigma$ stems from subcomplexes that are contractible, but not collapsible.
The most prominent such example is the $2$-dimensional \emph{dunce hat} which can be obtained from a single triangle by identifying 
its three boundary edges in a non-coherent way~\cite{Zeeman1963b}.

Crowley et al.~\cite{KCrowleyEbinKahnReyfmanWhiteXue2003pre} showed that the $7$-simplex with $8$~vertices contains in its $2$-skeleton an $8$-vertex 
triangulation of the dunce hat onto which it collapses.
(This result can easily be verified by running the random discrete Morse Step~(2) on the $7$-simplex, but not allowing
the triangles of an $8$-vertex triangulation of the dunce hat be used as free faces.)
While the dunce hat has triangulations with $8$ vertices~\cite{BenedettiLutz2013b},
every contractible complex with fewer vertices is collapsible~\cite{BagchiDatta2005}.

This leads us to our next experiment, where we compute random discrete Morse vectors for $d$-simplices, $7 \leq d \leq 25$; cf.\ Table~\ref{tbl:dsimplex}.
For instance, in dimension seven every one of the $10^{10}$ runs that we tried gave a perfect discrete Morse function, i.e., a collapsing sequence.
With increasing dimension that success rate drops slowly until $d=20$, where we get stuck with a discrete Morse vector which is not perfect in $1.3\%$ of all tries.
Going to even higher dimensions shows a rapid decline of the probability to find a perfect discrete Morse vector.
From this we conclude that dimension $25$ marks a \enquote{horizon} for Step~(2), even for a single simplex.
Note also that the implementation of the algorithm requires to store the entire Hasse diagram, which is memory expensive; 
e.g., the Hasse diagram of the $25$-simplex needs around $200$ GB of~RAM.
The data structure underlying the Hasse diagram is optimized for speed with respect to enumerating all faces of a simplicial complex (and the covering relations of the inclusion poset) from the facets.
That algorithm, which is linear in the size of the output, builds on a very general method of Ganter~\cite{Ganter:1984} for closure systems, and the implementation details in \polymake are discussed in \cite{HampeJoswigSchroeter:MEGA2017}.

\begin{table}[t]
\small\centering
\defaultaddspace=0.15em
\caption{Discrete Morse vectors for $10^9$ runs on the $8$-simplex.}\label{tbl:8simplex}
\begin{tabular*}{\linewidth}{@{\extracolsep{\fill}}l@{\hspace{1mm}}r@{}}
\\\toprule
 \addlinespace
 Discrete Morse vector  &     Count   \\\midrule
 (1 0 0 0 0 0 0 0 0)     &   999999988 \\
 (1 1 1 0 0 0 0 0 0)     &           4 \\
 (1 0 1 1 0 0 0 0 0)     &           7 \\
 (1 0 0 1 1 0 0 0 0)     &           1 \\
 \addlinespace
\bottomrule
\end{tabular*}
\end{table}

\begin{figure}[b]\centering
    \hspace{-2cm}
  \begin{minipage}{.6\textwidth}\centering
    \includegraphics[height=2.4cm]{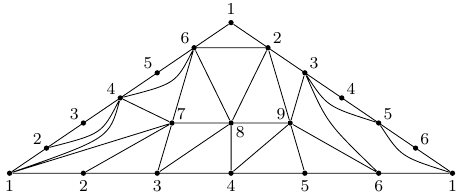}
  \vspace{0.5cm}
    
    \includegraphics[height=2.4cm]{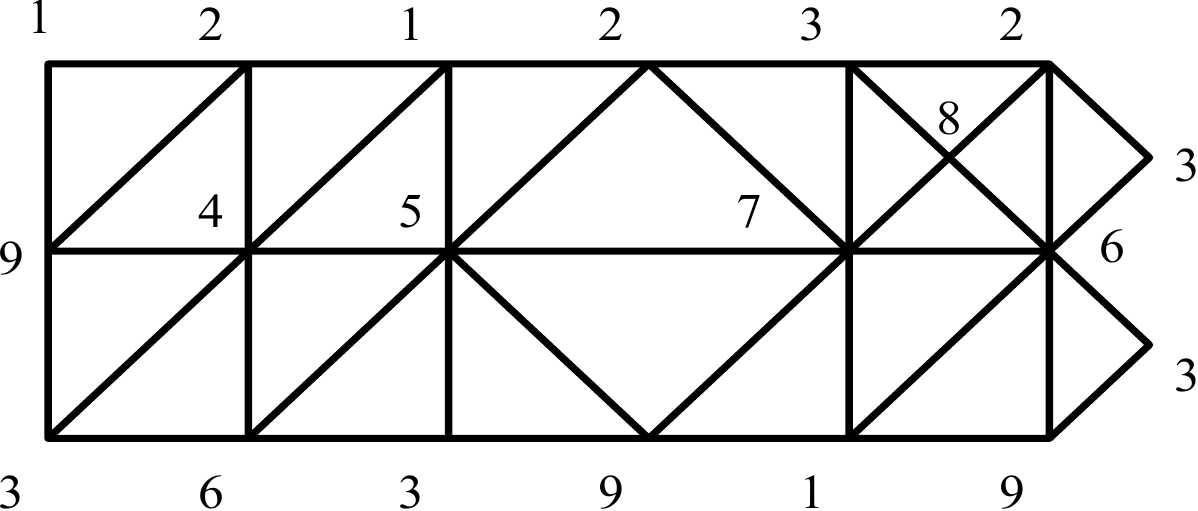}
  \end{minipage}
  \qquad
  \hspace{-1cm}
  \begin{minipage}{.3\textwidth}\centering
      \includegraphics[height=6cm]{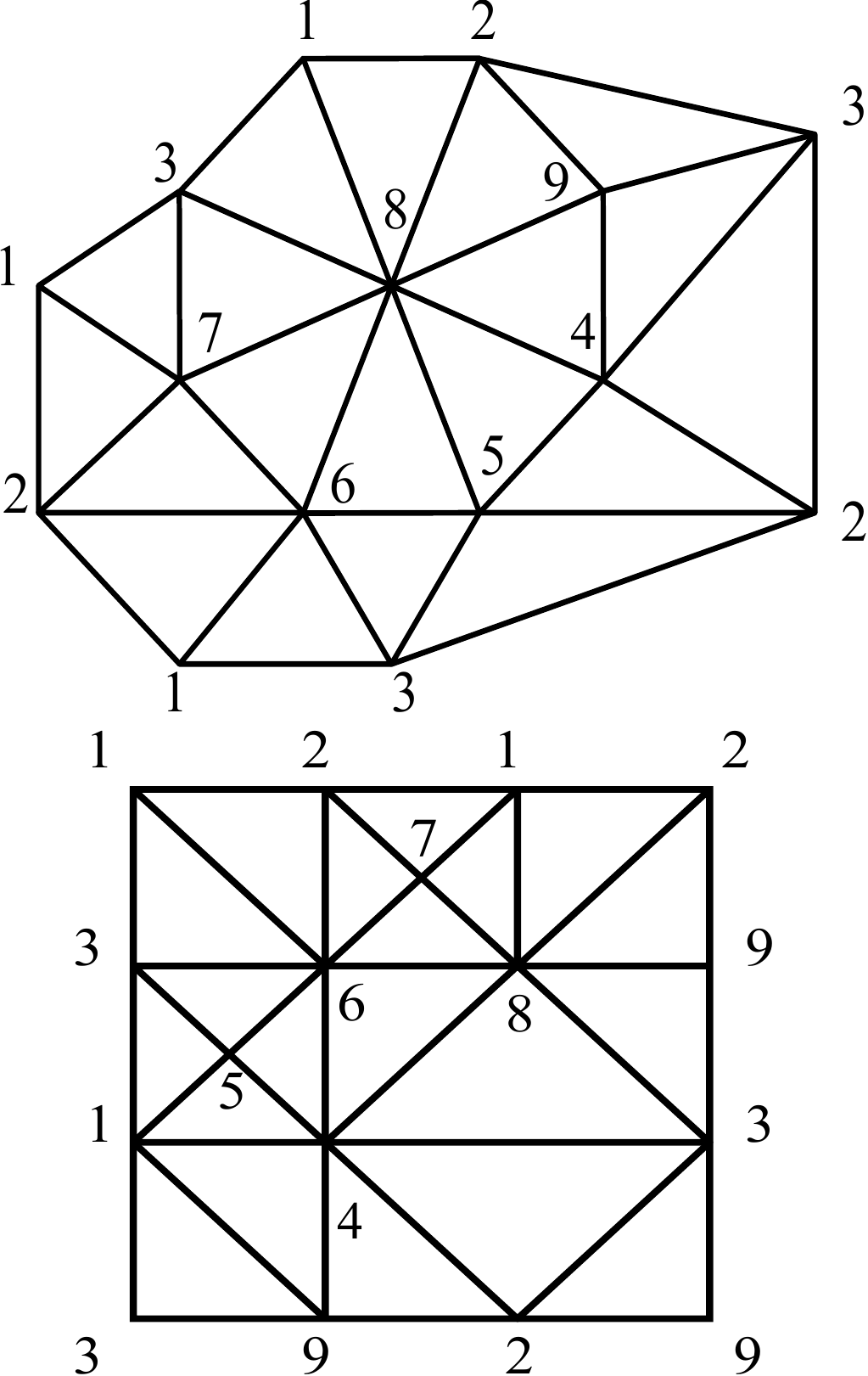}  
  \end{minipage}
  \caption{The simplicial complexes $\DunceHat$, $\Sawblade^2_a$, $\Sawblade^3_b$ and $\SawSquare$ (clockwise from top left).}
  \label{fig:sawblade-examples}
\end{figure}

It is instructive to look at the subcomplexes which arise from non-perfect discrete Morse vectors  in our experiments for a $d$-simplex.
For  $d=8$, we found four examples of $2$-dimensional contractible and non-collapsible complexes on nine vertices 
which we call $\DunceHat$, $\Sawblade^2_a$, $\Sawblade^3_b$, and $\SawSquare$; cf.\ Figure~\ref{fig:sawblade-examples} and (the second line of) Table~\ref{tbl:8simplex}.
The first one, $\DunceHat$, is a triangulation of the dunce hat.
The following concept is derived from scrutinizing the complexes in Figure~\ref{fig:sawblade-examples}.

\begin{figure}[b]
\begin{center}
\includegraphics[width=0.4\columnwidth]{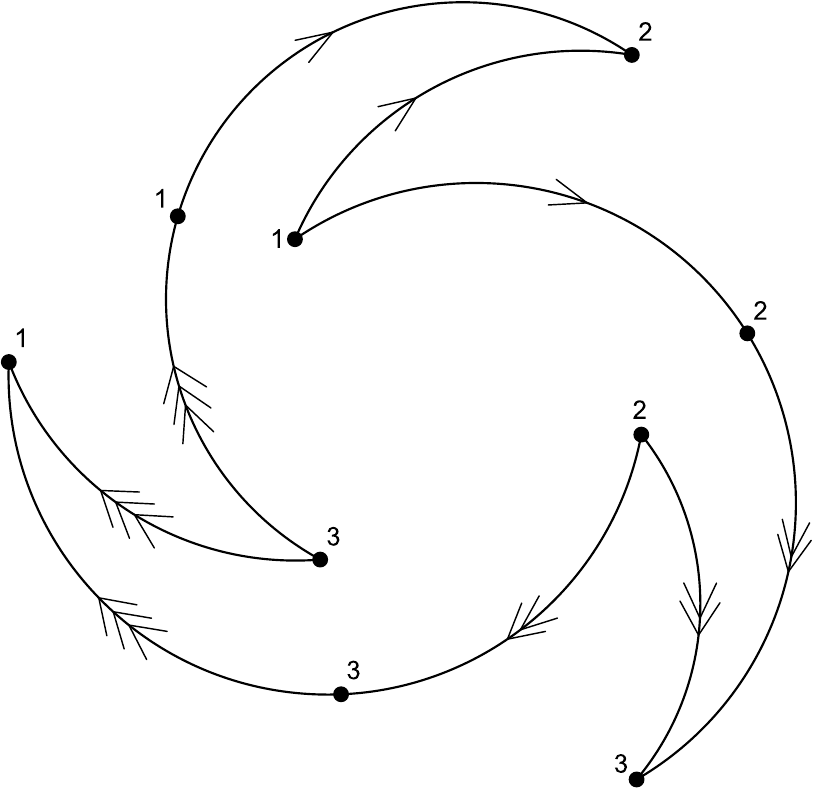}
\end{center}
\caption{The sawblade complex $\Sawblade^3$.} 
\label{fig:sawblade-cw}
\end{figure}

\begin{deff}\label{deff:saw_blade}
  The \emph{$k$-bladed saw blade complex $\Sawblade^k$} is the $2$-dimensional CW complex obtained from a polygonal disk with $3k$ edges $a_1,a_1^{-1},a_1, a_2,a_2^{-1},a_2,\dots, a_k,a_k^{-1},a_k$ by identifying $a_1,a_1^{-1},a_1$ as well $a_2,a_2^{-1},a_2$ and so on until $a_k,a_k^{-1},a_k$, for $k\geq 1$.
\end{deff}
In particular, for $k=1$ we obtain a triangle whose three edges are identified in the order $a_1,a_1^{-1},a_1$, i.e., $\Sawblade^1$ is the dunce hat.
More generally, $\Sawblade^k$ consists of $k$ vertices, $k$ edges, and a single disk; so the Euler characteristic equals one; cf.\ Figure~\ref{fig:sawblade-cw}, which explains the name.
We use notation like $\Sawblade^k_x$ and $\Sawblade^k_y$ to denote specific triangulations.

\begin{thm}
The following holds for $k$-bladed saw blade complexes $\Sawblade^k$:
\begin{itemize}
 \item[(i)] The dunce hat $\Sawblade^1$ can be triangulated with $8$ vertices.
 \item[(ii)] $\Sawblade^2$ can be triangulated with $9$ vertices.
 \item[(iii)] $\Sawblade^k$ can be triangulated with $3k$ vertices, for $k\geq 3$.
 \item[(iv)] Any triangulation of a saw blade complex is contractible, but non-collapsible.
\end{itemize}
\end{thm}

\emph{Proof:} We first prove (iii): For $k\geq 3$, assume that the identified boundary of the polygonal disk reads $1$--$2$--$1$--$2$--$3$--$2$--$3$--$4$--$3$--$4$--\dots--$k$--$1$--$k$--$1$; 
    see Figure~\ref{fig:sawblade-cw} in the case $k=3$.
    In the interior of the disk we place a cycle with $2k$ vertices and connect the cycle vertices with the boundary cycle vertices.
    More precisely, we connect every other cycle vertex to the beginning vertex of a blade and its two neighbors, and we connect the remaining cycle vertices to the two middle vertices of a blade.
    Finally, the interior $2k$-gon can be triangulated arbitrarily without additional vertices.
    
(ii)+(i): In the case of two blades, we start with $1$--$2$--$1$--$2$--$3$--$1$--$3$--$2$--$3$--$1$ as the identified boundary cycle.
    The extra vertex, say 3, is needed to avoid unwanted additional identifications.
    In the interior we then place a hexagon and connect its vertices similar to before.
    For an $8$-vertex triangulation of the dunce hat $\Sawblade^1$ see~\cite{BenedettiLutz2013b}.
    
(iv): The dunce hat $\Sawblade^1$ is contractible, and none of its triangulations is collapsible~\cite{Zeeman1963b}.
    For $k\geq 2$, the saw blade complex $\Sawblade^k$ has $k$ vertices that (by labeling appropriately) appear in order \mbox{$1$--$2$--$1$--$2$--$3$--$2$--$3$--$4$--$3$--$4$--\dots--$k$--$1$--$k$--$1$} on the boundary of the original $3k$-gon.
    We cut the $3k$-gon along an interior arc into two polygonal disks with identifications on the boundary, $1$--$2$--$1$--$2$--(interior arc) and the remainder $2$--$3$--$2$--$3$--$4$--$3$--$4$--\dots--$k$--$1$--$k$--$1$--(interior arc). 
    Both parts are contractible CW complexes (that retract to the paths $1$--$2$ and $2$--$3$--$4$--\dots--$k$--$1$) glued along a contractible subcomplex (the interior arc).
    We conclude that their union $\Sawblade^k$ is contractible; cf.\ Hachimori~\cite{Hachimori2008} for a similar decomposition of a contractible $2$-complex.

    Now we consider any triangulation $K$ of the saw blade complex $\Sawblade^k$, for $k\geq 1$.
    There is no free $2$-face as all edges in such a triangulation either have degree two or three. 
    It follows that $K$ is non-collapsible. \hfill $\Box$

\bigskip

Saw blade triangulations with different numbers of blades are combinatorially non-isomorphic.
These simplicial complexes and their quotients provide $2$-dimensional contractible, but non-collapsible simplicial complexes on which we can get stuck when trying to randomly collapse a simplex.
For instance, the simplicial complex $\SawSquare$ from Figure~\ref{fig:sawblade-examples} is obtained as a quotient from identifying two vertices in some triangulation of $\Sawblade^3$.
Similar examples and higher-dimensional analogs exist in abundance.

\begin{table}[t]
\small\centering
\defaultaddspace=0.15em
\caption{Discrete Morse vectors for $10^3$ runs on the $20$-simplex.}\label{tbl:20simplex}
\begin{tabular*}{\linewidth}{@{\extracolsep{\fill}}l@{\hspace{1mm}}r@{}}
\\\toprule
 \addlinespace
 Discrete Morse vector  &  Count   \\\midrule
 (1, 0, 0, 0, 0, 0, 0, 0, 0, 0, 0, 0, 0, 0, 0, 0, 0, 0, 0, 0, 0)                                  &  987  \\
 (1, 0, 0, 0, 6, 26, 59, 87, 61, 13, 0, 0, 0, 0, 0, 0, 0, 0, 0, 0, 0)                        &      1  \\
 (1, 0, 3, 30, 111, 158, 132, 82, 24, 0, 0, 0, 0, 0, 0, 0, 0, 0, 0, 0, 0)                &      1  \\
 (1, 0, 1, 8, 34, 80, 126, 155, 126, 61, 27, 10, 0, 0, 0, 0, 0, 0, 0, 0, 0)            &     1  \\
 (1, 0, 1, 14, 27, 24, 13, 3, 0, 0, 0, 0, 0, 0, 0, 0, 0, 0, 0, 0, 0)                          &     1  \\
 (1, 0, 1, 30, 117, 278, 409, 393, 213, 39, 0, 0, 0, 0, 0, 0, 0, 0, 0, 0, 0)          &     1  \\
 (1, 0, 2, 25, 110, 236, 305, 175, 19, 0, 0, 0, 0, 0, 0, 0, 0, 0, 0, 0, 0)              &     1  \\
 (1, 3, 5, 9, 34, 85, 134, 109, 33, 0, 0, 0, 0, 0, 0, 0, 0, 0, 0, 0, 0)                    &     1 \\
 (1, 0, 1, 19, 82, 150, 161, 90, 15, 0, 0, 0, 0, 0, 0, 0, 0, 0, 0, 0, 0)                  &     1 \\
 (1, 0, 3, 18, 51, 118, 196, 264, 207, 57, 0, 0, 0, 0, 0, 0, 0, 0, 0, 0, 0)            &     1 \\
 (1, 0, 1, 11, 107, 243, 366, 463, 450, 261, 54, 0, 0, 0, 0, 0, 0, 0, 0, 0, 0)      &     1 \\
 (1, 0, 1, 5, 30, 95, 160, 163, 124, 72, 27, 7, 0, 0, 0, 0, 0, 0, 0, 0, 0)              &     1 \\
 (1, 0, 6, 48, 182, 377, 657, 876, 801, 493, 170, 22, 0, 0, 0, 0, 0, 0, 0, 0, 0)  &     1 \\
 (1, 0, 0, 0, 0, 0, 0, 8, 14, 13, 14, 7, 0, 0, 0, 0, 0, 0, 0, 0, 0)                            &     1 \\
 \addlinespace
\bottomrule
\end{tabular*}
\end{table}

Tables~\ref{tbl:8simplex} and \ref{tbl:20simplex} give the actual discrete Morse vectors
found for the $8$-sim\-plex and the $20$-simplex, respectively. We observe 
that we can get stuck (i.e., run out of free faces) in different dimensions; see \cite{LofanoNewman2019} for an analysis of this phenomenon.

While in the case of the $8$-simplex at most two extra critical cells are picked up (see Table~\ref{tbl:8simplex}),
the discrete Morse vector $(1, 0, 6, 48, 182, 377, 657, 876, 801, 493, 170, 22, 0, 0, 0, 0, 0, 0, 0, 0, 0)$ 
for the $20$-simplex in Table~\ref{tbl:20simplex} contains $3,632$ extra critical cells.
Thus, in higher dimensions, not only do we get stuck with non-collapsible, contractible subcomplexes more often, but when we get stuck, the resulting discrete Morse vectors will also be larger. 
This may be seen as empirical evidence for the non-approximability of perfect Morse function; cf.\ \cite{BauerRathod2019}.

\subsection{Iterated barycentric subdivisions}
\label{sec:limit_morse_two}

As already seen in Section \ref{sec:runtimes}, it is sometimes rather easy to find perfect discrete Morse vectors, even in fairly large simplicial complexes, 
provided that the complexes are nicely structured; cf.\ also~\cite{BenedettiLutz2014}.
Adiprasito and Izmestiev~\cite{AdiprasitoIzmestiev} showed that sufficiently large iterated barycentric subdivisions of any PL sphere admit spherical discrete Morse functions.
Yet, the average number of critical cells for random discrete Morse vectors grows exponentially with the number of barycentric subdivisions~\cite{AdiprasitoBenedettiLutz2014pre}.
Here we try our sphere recognition heuristic on higher barycentric subdivisions of boundaries of simplices.

For the third barycentric subdivision $\sd^3\partial\Delta_4$ of the boundary of the $4$-simplex with face vector $(12600,81720,138240,69120)$ 
the perfect discrete Morse vector $(1,0,0,1)$ was found in $994$ out of $1000$ runs of the \texttt{random-lex-last} 
version~(cf.~\cite{AdiprasitoBenedettiLutz2014pre} and Section~\ref{subsec:randomMorse})
of the random discrete Morse search; see Table~\ref{tbl:discrete_morse_spectra}.
For $\sd^4\partial\Delta_4$ with face vector $(301680,1960560,3317760,1658880)$ the (same) perfect discrete Morse vector was found 
in only $844$ out of $1000$ runs. This suggests that the $4$th barycentric subdivision is still within the \enquote{horizon} for computations 
with the version \texttt{random-lex-last}, while the $5$th barycentric subdivision was too large to fit into the main memory of the machine we used for the experiments.
The \texttt{random-lex-first} strategy behaved slightly better than \texttt{random-lex-last}; 
the strategy \texttt{random-random} was always successful.

\begin{table}[t]
\small\centering
\defaultaddspace=0.15em
\caption{Discrete Morse vectors for iterated barycentric subdivisions of the $3$-sphere $\partial\Delta_4$.}\label{tbl:discrete_morse_spectra}
\begin{tabular*}{\linewidth}{@{}l@{\extracolsep{5mm}}r@{\extracolsep{\fill}}l@{\extracolsep{5mm}}r@{\extracolsep{\fill}}l@{\extracolsep{5mm}}r@{}}
\\\toprule
 \addlinespace
  \multicolumn{2}{@{}l@{}}{\texttt{random-random}}      &  \multicolumn{2}{@{}l@{}}{\texttt{random-lex-first}} &  \multicolumn{2}{@{}l@{}}{\texttt{random-lex-last}} \\ \midrule
 \addlinespace
\multicolumn{6}{@{}c@{}}
  {$\sd^3\partial\Delta_4$\,\,  with\, $f=(12600,81720,138240,69120)$} \\
  \midrule
  \addlinespace
 \textbf{(1,0,0,1)}:     &  1000 & \textbf{(1,0,0,1)}: &  999 & \textbf{(1,0,0,1)}: &  994 \\
                                                        &       & $(1,1,1,1)$:        &    1 & $(1,1,1,1)$:        &    6 \\
  \midrule
\multicolumn{6}{@{}c@{}}
  {$\sd^4\partial\Delta_4$\,\,  with\, $f=(301680,1960560,3317760,1658880)$} \\
  \midrule
  \addlinespace
 \textbf{(1,0,0,1)}:           &  100 & \textbf{(1,0,0,1)}: &  829 & \textbf{(1,0,0,1)}: &  844 \\
                         &    & $(1,1,1,1)$:        &  143 & $(1,1,1,1)$:        &  107 \\[-.5mm]
                         &    & $(1,2,2,1)$:        &   19 & $(1,2,2,1)$:        &   30 \\[-.5mm]
                         &    & $(2,3,2,1)$:        &    3 & $(1,3,3,1)$:        &    9 \\[-.5mm]
                         &    & $(2,5,4,1)$:        &    2 & $(1,4,4,1)$:        &    4 \\[-.5mm]
                         &    & $(1,3,3,1)$:        &    2 & $(2,5,4,1)$:        &    2 \\[-.5mm]
                         &    & $(1,4,4,1)$:        &    1 & $(1,5,5,1)$:        &    2 \\[-.5mm]
                         &    & $(1,5,5,1)$:        &    1 & $(2,3,2,1)$:        &    1 \\[-.5mm]
                         &    &                     &      & $(2,7,6,1)$:        &    1  \\[1mm]
 \addlinespace
\bottomrule
\end{tabular*}
\end{table}

\subsection{Other input}

Except for the Akbulut--Kirby spheres all the examples studied so far arise from easy to understand procedures.
Searching for entirely different triangulations of $\Sph^3$, we started out with $\partial\Delta_4$ with five vertices.
Then we added 525 vertices via random $0$-moves, followed by 50,000 random $1$-moves, followed by another $10^6$ rounds of random bistellar moves where we allowed both $1$- and $2$-moves.
This resulted in a \enquote{random} triangulation of the $3$-sphere with face vector $f=(530, 50474, 99888, 49944)$, which we fed into our heuristic.
The smallest discrete Morse vector found was $(1,2192,2192,1)$---far away from the perfect vector $(1,0,0,1)$.
This means that Step (2) fails on such input. Yet, applying bistellar moves again quickly gives back the initial $\partial\Delta_4$.
We also used \GAP to actually find a trivial presentation for the fundamental group of the example, which took 16 hours for the simplification.
It could be interesting to further investigate this or similar classes of examples.

\section*{Acknowledgments}
\noindent
We are indebted to Niko Witte for his original \polymake implementation of the bistellar simplification algorithm.
Further, we are grateful to Bruno Benedetti, Michael Kerber, Konstantin Mischaikow, Vidit Nanda, and five anonymous reviewers for their valuable comments.
In particular, we would like to thank one of the reviewers, who suggested to us to process the $4$-manifold census of \regina.
This lead to the results of Section~\ref{sec:4-manifolds}.

{M.~Joswig is supported by Einstein Foundation Berlin and Deutsche Forschungsgemeinschaft (EXC 2046: \enquote{MATH$^+$}, SFB-TRR 195: \enquote{Symbolic Tools in Mathematics and their Application}, GRK 2434: \enquote{Facets of Complexity}).
D.~Lofano is supported by the GRK 2434: ``Facets of Complexity''.
F.~H.~Lutz was partially supported by \textsc{VILLUM FONDEN} through the Experimental Mathematics Network and by the Danish National Research Foundation (DNRF) through the Center for Symmetry and Deformation. 
M.~Tsuruga was supported by the Berlin Mathematical School (BMS).

\bibliographystyle{amsplain}
\bibliography{recognition}


\end{document}